\documentclass[12pt]{article}
\usepackage{latexsym}
\usepackage{amsfonts}
\usepackage[all]{xy}
\usepackage{amssymb, mathrsfs, amsfonts, amsmath}
\usepackage{amsbsy}
\newcommand{\grad}{{{\rm{grad}\,}}}
\newcommand{\Hess}{{{\rm Hess}\,}}

\newcommand{\hessian}{{{\rm Hess}\,}}
\newcommand{\R}{\mathbb{R}}
\newtheorem{theorem}{\bf Theorem}[section]

\newtheorem{corollary}[theorem]{\bf Corollary}

\def\qed{\ifhmode\unskip\nobreak\fi\ifmmode\ifinner\else
\hskip5 pt \fi\fi\hbox{\hskip5 pt \vrule width4 pt height6 pt
depth1.5 pt \hskip 1pt }}

\begin{document}
\title{The mean curvature of cylindrically bounded submanifolds}
\author{L. J. Al\'{\i}as\thanks{Partially supported by MEC projects MTM2007-64504 and PCI2006-A7-0532, and Fundaci\'{o}n S\'{e}neca
project 04540/GERM/06, Spain. This research is a result of the activity developed within the framework of the Programme in Support of Excellence Groups of the Regi\'{o}n de Murcia, Spain, by Fundaci\'{o}n S\'{e}neca, Regional Agency
for Science and Technology (Regional Plan for Science and Technology 2007-2010).}
\and G. Pacelli Bessa\thanks{Partially supported by CNPq (Brazil), MEC project PCI2006-A7-0532 (Spain), and The Abdus Salam Int. Centre for Theoretical Physics-ICTP.}
\and  M. Dajczer\thanks{Partially supported by CNPq and Faperj (Brazil), and MEC project PCI2006-A7-0532 (Spain)}}
\date{Dedicated to Professor Manfredo P. do Carmo\\ on the occasion of his 80th birthday.}

\maketitle

\begin{abstract}
We give an estimate of the mean curvature  of a complete submanifold lying inside a closed cylinder $B(r)\times\R^{\ell}$ in a product Riemannian  manifold $N^{n-\ell}\times\R^{\ell}$. It follows that a complete  hypersurface of given
constant mean curvature lying inside a closed circular cylinder in Euclidean space cannot be proper if the circular base is of sufficiently small radius.  In particular, any possible counterexample to a conjecture of Calabi on
complete minimal hypersurfaces cannot be proper. As another application of our method, we derive a result about the stochastic incompleteness of submanifolds with sufficiently small mean curvature.

\

\end{abstract}

\section{Introduction}

\hspace*{0,6cm} The Calabi problem  in its original form, presented by Calabi \cite{calabi2} and promoted by Chern \cite{chern} about the same time, consisted on two conjectures  about Euclidean minimal hypersurfaces. The first
conjecture is that any complete minimal hypersurface of $\R^{n}$ must be unbounded. The second and more ambitious  conjecture asserted that any complete non-flat minimal hypersurface in $\R^{n}$ has unbounded projections in every
$(n-2)$-dimensional subspace.

Both conjectures turned out to be false for immersed surfaces in $\R^{3}$.  First Jorge and Xavier \cite{jorge-xavier2} exhibit a non-flat complete minimal surface lying between two parallel planes. Later on Nadirashvilli
\cite{Nadirashvili}  constructed a complete minimal surface inside a round ball in $\R^{3}$.

It was recently  shown by Colding and Minicozzi \cite{colding-minicozzi} that both conjectures hold for embedded minimal surfaces. Their work involves the close relation between the Calabi conjectures and properness. Recall that an
immersed submanifold in Euclidean space is \emph{proper} if the pre-image of any compact subset of $\R^n$ is compact. It is a consequence of their general result that a complete embedded minimal disk in $\R^{3}$ must be proper.

The immersed counterexamples to Calabi's conjectures discussed above  are not proper.  The example of Nadirashvilli cannot be proper since from the definition a proper submanifold must be unbounded. The same conclusion hold for the
other example but now the argument is not so easy, one has to use the strong half-space theorem due to Hoffman and Meeks \cite{H-M}.

The strong half-space theorem does not hold in $\R^n$ for $n\geq 4$. In fact, the higher dimensional catenoids are between parallel hyperplanes. Hence, it is natural to ask if  any possible higher dimensional counterexample to
Calabi's second conjecture must be  non-proper. In the special case of  minimal immersion, it follows from the corollary of our main result that a complete hypersurface of $\R^n$, $n\geq 3$, with bounded projection in a two
dimensional subspace cannot be  proper (see Corollary \ref{thm1} below).

As an application of our method, we  generalize the results by Markvorsen \cite{markvorsen} and Bessa and Montenegro \cite{bessa-montenegro1} about stochastic incompleteness of minimal submanifolds to submanifolds of bounded mean
curvature. In this respect, let us recall that a Riemannian manifold $M$ is said to be \emph{stochastically complete} if for some (and therefore, for any) $(x,t)\in M\times(0,+\infty)$ it holds that $\int_Mp(x,y,t)dy=1$, where
$p(x,y,t)$ is the heat kernel of the Laplacian operator. Otherwise, the manifold $M$ is said to be \emph{stochastically incomplete} (for further details about this see, for instance, \cite{Gr} or \cite{pigola-rigoli-setti}).

An interesting problem in submanifold geometry is to understand stochastic completeness/incompleteness of submanifolds in terms of their extrinsic geometry. In \cite{markvorsen} Markvorsen derived a mean time exit comparison theorem
which implies that any bounded complete minimal submanifold of a Hadamard manifold $N$ with sectional curvature $K_{N}\leq b\leq 0$ is stochastically incomplete. Recently,  Bessa and Montenegro \cite{bessa-montenegro1} considered
minimal submanifolds of product spaces $N\times \mathbb{R}$, where $N$ is a Hadamard manifold with $K_{N}\leq b\leq 0$, and proved a version of Markvorsen's result in this setting. In particular, they showed that complete
cylindrically bounded  minimal submanifolds of $N\times \mathbb{R}$ are stochastically incomplete. Here we extend these results to complete submanifolds with sufficiently small mean curvature lying inside a closed cylinder
$B(r)\times\R^{\ell}$ in a product Riemannian manifold $N^{n-\ell}\times\R^{\ell}$.

\section{The results}
\hspace*{0,6cm}

Part (a) of Theorem \ref{thm3} below extends the main results given in  \cite{bessa-montenegro} for compact hypersurfaces.
Part (b)  generalizes stochastic incompleteness results of \cite{bessa-montenegro1} and \cite{markvorsen} for minimal
submanifolds.

In the following we denote
$$
C_b(t)
=\left\{\begin{array}{lll}
\sqrt{b}\cot(\sqrt{b}\, t) & \mathrm{if} & b >0,\;\; t<\pi/2\sqrt{b},\\
1/t & \mathrm{if} & b =0,\\
\sqrt{-b}\coth(\sqrt{-b}\, t) & \mathrm{if}  & b <0.
\end{array}\right.
$$

\begin{theorem}\label{thm3}
Let $\varphi\colon M^m\to N^{n-\ell}\times\R^{\ell}$ be an isometric immersion of a complete Riemannian manifold $M$ of
dimension $m\geq \ell+1$. Let $B_N(r)$ be the geodesic ball of $N^{n-\ell}$ centered at $p$ with radius $r$. Given $q\in M$, assume that the radial sectional
curvatures $K_N^{\mathrm{rad}}$ along the radial geodesics issuing from
$p=\pi_N(\varphi(q))\in N^{n-\ell}$ are bounded
as $K_N^{\mathrm{rad}}\leq b$ in $B_N(r)$.  Suppose that
$$
\varphi (M)\subset B_N(r)\times \R^\ell
$$
for $r<\min\{{\rm inj}_N(p), \pi/2\sqrt{b}\}$, where we replace $\pi/2\sqrt{b}$ by $+\infty$ if $b\leq 0$.
\begin{enumerate}
\item[(a)] If $\varphi\colon M^m\to N^{n-\ell}\times\R^{\ell}$ is proper, then
\begin{equation}
\sup_{M}| H| \geq \frac{(m-\ell)}{m}C_{b}(r).
\end{equation}
\item[(b)] If
\begin{equation}
\sup_{M}| H|<\frac{(m-\ell)}{m}C_{b}(r),
\end{equation}
then $M$ is stochastically incomplete.
\end{enumerate}
\end{theorem}

For Euclidean hypersurfaces we have the following consequence.
\begin{corollary}\label{thm1}
Let $\varphi\colon M^{n-1}\to\R^{n}$ be a complete hypersurface with mean curvature $H$.  If $\varphi(M)\subset B_{\R^2}(r)\times \R^{n-2}$ and $\sup_M|H| < 1/(n-1)r$, then $\varphi$ cannot be proper.
\end{corollary}
Observe that the assumption on the bound of the mean curvature cannot be weakened since $1/(n-1)r$ is the mean curvature of the cylinder $\mathbb{S}^1(r)\times \R^{n-2}$.

We point out that Mart\'{\i}n and Morales \cite{pacomartin} constructed  examples of complete  minimal surfaces properly immersed in the interior of a cylinder $B_{\R^2}(r)\times\R$. By the above
result these surfaces cannot be proper in $\R^3$.

\section{The proofs}
Let $\varphi\colon M^m \to N^n$ be an isometric immersion between Riemannian manifolds. Given a function $g\in C^\infty(N)$ we set $f=g\circ\varphi\in C^\infty(M)$.
Since
$$
\langle\grad^Mf,X\rangle=\langle\grad^Ng, X\rangle
$$
for every vector field $X\in TM$, we obtain
$$
\grad^Ng=\grad^Mf +(\grad^Ng)^{\perp}
$$
according to the decomposition $TN=TM\oplus T^\perp M$. An easy computation using the Gauss formula gives the well-known relation
(see e.g. \cite{jorge-koutrofiotis})
\begin{equation}\label{eqBF2}
\Hess_Mf(X,Y)= \Hess_Ng(X,Y) +\langle\grad^Ng,\alpha(X,Y)\rangle
\end{equation}
for all vector fields $X,Y\in TM$, where $\alpha$ stands for the second fundamental form of
$\varphi$. In particular, taking traces  with respect to an orthonormal frame
$\{ e_{1},\ldots, e_{m}\}$ in $TM$ yields
\begin{equation}\label{eqBF3}
\Delta_Mf =\sum_{i=1}^{m}\Hess_Ng(e_i,e_i)+ \langle\grad^Ng,\stackrel{\to}{H}\rangle.
\end{equation}
where $\stackrel{\to}{H}=\sum_{i=1}^{m}\alpha(e_i,e_i)$.
\vspace{1,5ex}

The first main ingredient of our proofs is the Hessian comparison theorem.

\begin{theorem}\label{thm2}
Let $M^m$ be a  Riemannian manifold  and
$x_0,x_1 \in M$ be such that there is a minimizing unit speed geodesic $\gamma$  joining $x_{0}$ and $x_{1}$ and let $\rho(x)=\mathrm{dist}(x_0,x)$
be the  distance function  to $x_{0}$. Let $K_{\gamma}\leq b$ be the radial
sectional curvatures of $M$ along $\gamma$.
If $b>0$ assume $\rho(x_{1})<\pi/2\sqrt{b}$. Then, we have $\Hess\rho(x) (\gamma',\gamma')=0$ and
\begin{equation}\label{thmHess}
\Hess\rho(x)(X,X)\geq C_{b}(\rho(x))\Vert
X\Vert^2
\end{equation}
where $X\in T_{x}M$ is perpendicular to $\gamma'(\rho(x))$.
\end{theorem}

The second main ingredient is the version proved by Pigola-Rigoli-Setti \cite[Theorem 1.9]{pigola-rigoli-setti} of the
Omori-Yau maximum principle.

\begin{theorem}\label{eqThmMP}
Let $M^m$ be a  Riemannian manifold and assume that there exists a non-negative $C^2$-function $\psi$ satisfying the following requirements:
$$
\psi(x)\rightarrow +\infty \;\;\; {\rm as} \;\;\; x\rightarrow \infty $$
$$
\exists\, A>0\;\;\; {\rm such\,\, that}\;\;\; |\grad\psi|\leq  A\sqrt{\psi}\;\;\;\hspace*{6,5ex}{\rm off\,\, a\,\, compact\,\, set}
$$
$$
\exists\, B>0 \;\;\; {\rm such\,\, that} \;\;\; \Delta\psi \leq  B\sqrt{\psi G(\sqrt{\psi})}\;\;\;{\rm off \,\, a \,\, compact\,\, set}
$$
where $G$ is a smooth function on $[0, +\infty)$ satisfying:
\begin{equation}\begin{array}{llll}
(i)\,\, G(0)>0, &  & & (ii) \,\, G'(t)\geq 0 \,\,{\rm on}\,\,[0, +\infty),\vspace{1ex}\\
(iii)\,\,1/ \sqrt{G(t)}\not\in L^{1}(0,+\infty), &  & &
(iv) \,\, \limsup_{t\to +\infty}\displaystyle{\frac{tG(\sqrt{t})}{G(t)}}<+\infty.
\end{array}\label{eqG}
\end{equation}
Then, given a function $u\in C^2(M)$ with  $u^{\ast}=\sup_{M}u <+\infty$
there exists a sequence  $\{x_{k}\}_{k\in\mathbb{N}}\subset M^m$ such that
$$
\begin{array}{lllll}u(x_{k})>u^{\ast}-1/k;& & | \grad u  |(x_{k}) <1/k; && \Delta u (x_{k})<1/k. \end{array}
$$
\end{theorem}

Observe that a function $G$ satisfying the above conditions is
\begin{equation}
\label{G}
G(t)=(t+2)^2(\log(t+2))^2.
\end{equation}
Now we are ready to prove  Theorem \ref{thm3}.\vspace{1,5ex}

\noindent {\it Proof of Theorem \ref{thm3}:} Define $\sigma:N^{n-\ell}\times\mathbb{R}^\ell\rightarrow [0,+\infty)$
by
\[
\sigma(z,y)=\rho_{\mathbb{R}^\ell}(y),
\]
where $\rho_{\mathbb{R}^\ell}(y)=\|y\|_{\mathbb{R}^\ell}$ is the distance function to the origin in $\mathbb{R}^\ell$.
Since $\varphi $ is proper and $\varphi(M)\subset B_N(r)\times\mathbb{R}^\ell$, then the function
$\psi(x)= \sigma\circ \varphi (x)$ satisfies $\psi(x)\to \infty$ as $\rho_{M}(x)={\rm dist}_{M}(q,x)\to +\infty$.
Off a compact set, we now have
$$
|\grad^{M}\psi(x)| \leq
|\grad^{N\times \R^\ell}\sigma(\varphi(x))|
 = | \grad^{\R^\ell}\rho_{\mathbb{R}^\ell}|=1\leq\sqrt{\psi(x)}.
$$
To compute $\Delta_{M}\psi$ we start with bases  $\{\partial/\partial\rho_N,\partial/\partial\theta_{2},\ldots,\partial/\partial\theta_{n-\ell}\}$ of $TN$ and
$\{\partial/\partial\rho_{\mathbb{R}^\ell},\partial/\partial\gamma_{2},\ldots,\partial/\partial\gamma_{\ell} \}$ of $T\R^\ell$ (polar coordinates) orthonormal at $x\in M$. Then, we choose  an orthonormal basis $\{e_1,\ldots,e_m\}$
for $T_{x}M$ as follows
$$
e_i = \alpha_i \frac{\partial}{\partial \rho_N}
+ \displaystyle \sum_{j=2}^{n-\ell} a_{ij}
\frac{\partial }{\partial \theta_{j}}+ \beta_i\frac{\partial}{\partial\rho_{\mathbb{R}^\ell}}
+ \displaystyle\sum_{t=2}^{\ell} b_{it}
\frac{\partial}{\partial \gamma_t}\cdot
$$
Hence, we have
$$
\Hess_{N\times \R^\ell}\, \sigma(\varphi(x))(e_i,e_i)
=\Hess_{\R^{\ell}}\,\rho_{\mathbb{R}^\ell}(\pi_{\R^\ell}e_i,\pi_{\R^\ell}e_i)
=\frac{1}{\sigma(\varphi(x))}\sum_{t=2}^\ell b^2_{it}\leq\frac{1}{\psi(x)},
$$
where $\pi_{\R^\ell}$ denotes the orthogonal projection onto $T\R^\ell$. Here, we are using
$$
|e_i|=1=\alpha_i^2+ \sum_{j=2}^{n-\ell}a_{ij}^2+\beta_i^2+\sum_{t=2}^\ell b_{it}^2
$$
that yields $\sum_{t=2}^\ell b^2_{it}\leq 1$.

Since $\psi(x)\to \infty$ as $\rho_{M}(x)={\rm dist}_{M}(q,x)\to +\infty$, off a compact set we may assume that
$$
|\hspace{-1mm}\stackrel{\to}{H}\hspace{-1mm}| (x)
=m| H|(x)\leq\sqrt{\psi(x) G(\sqrt{\psi (x)})}
$$
where $G(t)$ is given by (\ref{G}). Otherwise, $\sup_M|H|=+\infty$ and there is nothing to prove. Besides, off
a compact set we also have that
$$
\frac{1}{\psi(x)}\leq \sqrt{\psi(x) G(\sqrt{\psi (x)})}.
$$
Hence,  from (\ref{eqBF3}) we have off a compact set that
\begin{eqnarray}
\Delta_{M}\psi (x)\!\!&=&\!\! \displaystyle \sum_{i=1}^{m}\Hess_{N\times \R^\ell}\, \sigma(\varphi(x))(e_i,e_i)+
\langle\grad^{N\times \R^\ell} \sigma(\varphi (x)), \stackrel{\to}{H}(x)\rangle\nonumber \\
\!\!&\leq &\!\!  \frac{m}{\psi(x)}+ m|H| (x) \nonumber \\
\!\!&\leq &\!\! (m+1)\sqrt{\psi(x)G(\sqrt{\psi(x)})}.\nonumber
\end{eqnarray}
Therefore, by Theorem \ref{eqThmMP} the Omori-Yau maximum principle holds on $M$.

Define  $\rho\colon N^{n-\ell}\times \R^\ell\to \R$ by
$$
\rho (z, y)=\rho_N(z)={\rm dist}_N(p,z)
$$
and
$u\colon M^m\to\R$ by
$$
u(x)=\rho\circ\varphi(x).
$$
Since $\varphi (M)\subset B_N(r)\times\R^\ell$, we have that $u^{\ast}=\sup_{M}u\leq r<\infty$,
Therefore, by the maximum principle there is a sequence $\{x_{k}\}_{k\in\mathbb{N}}\subset M^m$  such that
$$\begin{array}{lllll}
u(x_{k})>u^{\ast}-1/k;
&& |\grad u|(x_{k}) <1/k; && \Delta u (x_{k})<1/k.
\end{array}
$$
Hence, we have
\begin{equation}\label{first}
\frac{1}{k}> \Delta u (x_k)=\sum_{i=1}^{m}\Hess_{N\times \R^\ell}\rho (\varphi(x_k))(e_i,e_i)+
\langle \grad^{N\times \R^\ell}\rho(\varphi(x_k)),\stackrel{\to}{H}\!(x_k)\rangle
\end{equation}
where $\{e_{1},\ldots,e_{m}\}$ is an orthonormal basis for $T_{x_{k}}M$.  Start with an
orthonormal basis  $\{\partial/\partial \rho_N,\partial/\partial \theta_{2},\ldots,
\partial/\partial \theta_{n-\ell}\}$ for $TN$ and  standard coordinates $\{y_1,\ldots y_{\ell}\}$ for $\R^{\ell}$.
Then, choose  an orthonormal  basis for $T_{x_{k}}M$ as
follows
$$
e_i = \alpha_i \frac{\partial}{\partial \rho_N}
+ \displaystyle \sum_{j=2}^{n-\ell}a_{ij}\frac{\partial}{\partial \theta_j}
+\displaystyle\sum_{t=1}^{\ell}c_{it}\frac{\partial}{\partial y_t}\cdot
$$
Using Theorem \ref{thm2}, a straightforward computation yields
\begin{eqnarray}
\Hess_{N\times \R^\ell}\rho (\varphi (x_{k}))(e_i,e_i)
\!\!&=&\!\!\Hess_N\rho_N (z (x_k))(\pi_{TN}e_i,\pi_{TN}e_i)\nonumber \\
\!\!&=&\!\! \sum_{j=2}^{n-\ell}a_{ij}^2\Hess_N\rho_N(z(x_k))
(\partial/\partial \theta_j,
\partial/\partial \theta_{j})\nonumber\\
\!\!&\geq &\!\!\sum_{j=2}^{n-\ell}a_{ij}^2C_b(r)\\
\!\!&=&\!\!
\left(1-\alpha_i^2-\sum_{t=1}^\ell c_{it}^2\right)C_b(r) \nonumber
\end{eqnarray}
since
$$
|e_i| = 1=\alpha_i^2+ \sum_{j=2}^{n-\ell}a_{ij}^2+ \sum_{t=1}^\ell c_{it}^2,
$$
where  $\pi_{TN}$ denotes the orthogonal projection onto $TN$. Therefore,
\begin{equation}\label{second}
\sum_{i=1}^m\Hess_{N\times \R^\ell}\rho(\varphi (x_k))(e_i,e_i)\geq \Big(m-\sum_{i}\alpha_i^2
-\sum_{i,t}c_{it}^2\Big)C_b(r).
\end{equation}
At $x_{k}$, we have
$$
\grad^{N\times\mathbb{R}^\ell}\rho(\varphi(x_k))=\grad u(x_k) +(\grad^{N\times\mathbb{R}^\ell}\rho(\varphi(x_k)))^{\perp}
$$
and hence
\begin{equation}\label{third}
|\grad u|^2(x_k)=\sum_{i=1}^m\langle\frac{\partial}{\partial \rho_N},e_i\rangle=\sum_i\alpha_i^2<1/k^2.
\end{equation}
Taking into account $|\grad^{N\times \R^\ell}\rho|=|\grad^{N}\rho_N|=1$, from (\ref{first}) and (\ref{second}) we obtain
$$
\frac{1}{k}>  \Big(m-\sum_i\alpha_i^2
-\sum_{i,t}c_{it}^2\Big)C_b(r) -m\sup_M|H|.
$$
It follows using (\ref{third}) that
\begin{equation}\label{forth}
\frac{1}{k}+ \frac{C_{b}(r)}{k^2}
+ m\sup_{M}|H|\geq\Big(m-\sum_{i,t}c^2_{it}\Big)C_b(r).
\end{equation}
Observe now that
\[
\sum_{i,t}c^2_{it}=\sum_{t=1}^\ell\sum_{i=1}^mc_{it}^2=\sum_{t=1}^\ell|\grad (y_t\circ\varphi)|^2\leq\ell,
\]
since $|\grad (y_t\circ\varphi)|^2\leq|\grad^{\mathbb{R}^\ell}y_t|^2=1$. Thus,
\[
m-\sum_{i,t}c^2_{it}\geq (m-\ell)
\]
and  we have letting $k\to +\infty$ in (\ref{forth}) that
$$
m \sup_{M}|H| \geq (m-\ell)C_{b}(r).
$$
This concludes the proof of the first part of Theorem \ref{thm3}.\vspace{1,5ex}

For the proof of the second part, we make use of the following
characterization of stochastic completeness given in \cite{PRS} (see \cite[Theorem 3.1]{pigola-rigoli-setti}):
A Riemannian manifold $M$ is stochastically complete if and only if for every $u\in C^{2}(M)$ with
$u^{\ast}=\sup u<\infty$ there exists a sequence $\{x_{k}\}$ such that $u(x_{k})> u^{\ast}-1/k$ and
$\Delta u (x_k)<1/k$ for every $k\ge 1$.\vspace{1,5ex}

Suppose that $M$ is stochastically complete. Define $g\colon N^{n-\ell}\times\R^{\ell}\to\R$ by
\[
g(z,y)=\hat{g}(z)=\phi_b(\rho_N(z))
\]
where
$$
\phi_b(t)=\left\{ \begin{array}{lll}
1-\cos(\sqrt{b}\, t) & \mathrm{if} & b>0, \;\; t<\pi/2\sqrt{b},\\
t^{2} & \mathrm{if} &  b=0,\\
\cosh(\sqrt{-b}\, t) & \mathrm{if} & b<0.
\end{array}\right.
$$
Then  $f=g\circ\varphi$ is a smooth bounded function on $M$. Thus there exists a sequence of points $\{x_k\}$ in $M$ such
that
$$
f(x_k)>f^*-1/k\;\;\;\mbox{and}\;\;\; \Delta f(x_k)<1/k
$$
for $k\ge 1$, where $f^*=\sup_Mf\leq\phi_b(r)<\infty$. Similar as before, we have
\begin{eqnarray}
\!\Hess_{N\times \R^\ell}g (\varphi (x_{k}))(e_i,e_i)
\!\!&=&\!\!\Hess_N \hat{g}(z(x_k))(\pi_{TN}e_i,\pi_{TN}e_i)\nonumber \\
\!\!&=&\!\! \phi_b''(r_k)\alpha_{i}^{2}+\phi'_b(r_k)\sum_{j=2}^{n-\ell}a_{ij}^2\Hess_N\rho_N(z(x_k))
(\partial/\partial \theta_j,\partial/\partial \theta_{j})\nonumber\\
\!\!& \geq &\!\! \phi_b''(r_k)\alpha_{i}^{2}+\phi'_b(r_k)C_b(r_k)\sum_{j=2}^{n-\ell}a_{ij}^2\nonumber\\
\!\!&=&\!\! \phi_b''(r_k)\alpha_{i}^{2}+\phi'_b(r_k)C_b(r_k)\left( 1-\alpha_{i}^{2}-\sum_{t=1}^\ell c_{it}^{2}\right)\nonumber\\
\!\!&=&\!\! \phi'_b(r_k)C_b(r_k)\left( 1-\sum_{t=1}^\ell c_{it}^{2}\right)\nonumber
\end{eqnarray}
since $\phi''_b(t)-\phi'_b(t)C_b(t)=0$. Here, we are writing $r_k=\rho_N(z(x_k))$. Therefore,
\begin{eqnarray}\frac{1}{k}> \Delta f (x_{k})\!\!&=&\!\!
\sum_{i=1}^{m}\hessian_{N\times\mathbb{R}^\ell}g (e_{i},e_{i})+ \langle \grad^{N\times\mathbb{R}^\ell} g, \stackrel{\to}{H}\rangle \nonumber\\
\!\!&\geq &\!\! \phi'_b(r_k)C_b(r_k)\left(m-\sum_{i,t}c_{it}^{2}\right)+\phi'_b(r_k)\langle \grad^{N\times\mathbb{R}^\ell}\rho_N, \stackrel{\to}{H}\rangle \nonumber \\
\!\!&\geq &\!\!\phi'_b(r_k)\left((m-\ell)C_b(r_k)-m\sup |H|\right).\nonumber
\end{eqnarray}
Finally, since $\lim_{k\to \infty}\phi'_b(r_{k})>0$, letting $k\to\infty$ we have
$$
\sup |H|\geq\frac{(m-\ell)}{m}\,C_b(r).
$$

\noindent{\it Proof of Corollary \ref{thm1}:} If  $\varphi$ is proper in $\R^{n}$, from part (a) of Theorem \ref{thm3} we would have $|H|\geq 1/(n-1)r$, and that is a contradiction.\qed

\

\textit{Note added in proof.} After submission of this paper, we were informed by Rosenberg that he and Sa Earp proved in \cite[Corollary 4.1.1 and Remark  4.3.3]{RSa} that a complete real analytic hypersurface $M$ properly immersed
into $\mathbb{R}^{n}$ which is inside a generalized rotational Delaunay hypersurface $\mathcal{D}$ and has mean curvature satisfying $|H|\leq \mathcal{H}_\mathcal{D}$ must be $M=\mathcal{D}$. Here $\mathcal{H}_\mathcal{D}$ denotes the
constant mean curvature of the generalized rotational Delaunay hypersurface $\mathcal{D}$. Although not stated in \cite{RSa}, it follows from this that a complete minimal hypersurface in $\mathbb{R}^{n}$, $n\geq 3$, with bounded
projection in an $(n-1)$-dimensional subspace cannot be proper.

{\renewcommand{\baselinestretch}{1}
\hspace*{-20ex}\begin{tabbing}
\indent \=Luis J. Alias \\
\>Departamento de Matematicas \\
\>Universidad de Murcia \\
\>  Campus de Espinardo E-30100 -- Spain\\
\> ljalias@um.es
\end{tabbing}}
\vspace*{-4ex}

{\renewcommand{\baselinestretch}{1} \hspace*{-20ex}\begin{tabbing}
\indent \= Gregorio Pacelli Bessa\\
\> UFC - Departamento de Matematica \\
\> Bloco 914 -- Campus do Pici\\
\> 60455-760 -- Fortaleza -- Ceara -- Brazil\\
\> bessa@mat.ufc.br
\end{tabbing}}
\vspace*{-4ex}

{\renewcommand{\baselinestretch}{1} \hspace*{-20ex}\begin{tabbing}
\indent \= Marcos Dajczer\\
\> IMPA \\
\> Estrada Dona Castorina, 110\\
\> 22460-320 -- Rio de Janeiro -- Brazil\\
\> marcos@impa.br\\
\end{tabbing}}


\begin{thebibliography}{abcd}

\bibitem{bessa-montenegro} G. P. Bessa and J. F. Montenegro,
{\em On compact H-Hypersurfaces of $N\times\R$.} Geom. Dedicata. {\bf 127} (2007), 1--5.

\bibitem{bessa-montenegro1}G. P. Bessa \and J. Fabio Montenegro,
{\em Mean time exit and isoperimetric inequalities for minimal submanifolds of $N\times \mathbb{R}$.} To appear in Bull. London Math. Soc.
Available at http://arxiv.org/pdf/0709.1331

\bibitem{calabi2}E. Calabi, {\em Problems in Differential Geometry} (S. Kobayashi and J. Eells, Jr., eds.) Proc. of the United States-Japan Seminar in Differential Geometry, Kyoto, Japan, 1965, Nippon Hyoronsha Co. Ltd., Tokyo (1966) 170.

\bibitem{chern} S. S. Chern, {\em The Geometry of G-structures.} Bull. Amer. Math. Soc. {\bf 72} (1966), 167--219.

\bibitem{colding-minicozzi} T. Colding and W. Minicozzi II, {\em The Calabi-Yau conjectures for embedded surfaces.} Annals of Math. \textbf{161} (2005) 727--758.


\bibitem{RSa} R. Sa Earp and H. Rosenberg, {\em Some remarks on surfaces of prescribed mean curvature}.
Differential geometry, 123--148, Pitman Monogr. Surveys Pure Appl. Math., 52, Longman Sci. Tech., Harlow, 1991.


\bibitem{Gr} A. Grigor'yan, \textit{Analytic and geometric background of recurrence and non-explosion of the Brownian motion on Riemannian manifolds},
Bull. Amer. Math. Soc. (N.S.) {\bf 36 (1999)}, 135--249.

\bibitem{H-M} D. Hoffman and W. Meeks, {\em The Strong Half-space Theorem  for minimal surfaces.} Invent. Math. {\bf 101}  (1990) 373--377.

\bibitem{jorge-koutrofiotis}L.  Jorge  and  D. Koutrofiotis,
{\em An estimate for the curvature of
bounded submanifolds.} Amer. J. Math.  {\bf 103} (1980) 711--725.

\bibitem{jorge-xavier2} L. Jorge and F. Xavier, {\em A complete minimal surface in $\R^{3}$ between two parallel planes.} Ann. of Math. {\bf 112} (1980) 203--206.

\bibitem{markvorsen} S. Markvorsen, {\em On the mean exit time from a
submanifol.} J. Differential Geom. \textbf{29} (1989), 1-8.

\bibitem{pacomartin}F. Mart\'{\i}n and S. Morales, {\em A complete bounded minimal cylinder in $\mathbb R\sp 3$. } Michigan Math. J.  {\bf 47}  (2000), 499--514.

\bibitem{Nadirashvili}N.  Nadirashvili, {\em Hadamard's and
Calabi-Yau's conjectures on negatively curved and minimal surfaces.}
Invent. Math. {\bf 126} (1996), 457--465.

\bibitem{PRS} S. Pigola, M. Rigoli and A. Setti,
{\em A remark on the maximum principle and stochastic completeness.} Proc. Amer. Math. Soc. {\bf 131} (2003), 1283--1288.

\bibitem{pigola-rigoli-setti} S. Pigola, M. Rigoli and A. Setti,
{\em Maximum Principle on Riemannian Manifolds ans Applications.} Memoirs Amer. Math. Soc. {\bf 822} (2005).
\end{thebibliography}
\end{document}